\documentclass[11pt,mathsrc]{amsproc}
\usepackage[dvips]{graphicx,color}
\usepackage{amssymb,amscd,latexsym,epsfig,euscript}
\usepackage[arrow,matrix,line]{xy}

\DeclareFontFamily{OT1}{rsfs}{}
\DeclareFontShape{OT1}{rsfs}{n}{it}{<-> rsfs10}{}
\DeclareMathAlphabet{\curly}{OT1}{rsfs}{n}{it}

\newcommand\C{\mathbb C}
\newcommand\R{\mathbb R}
\newcommand\Z{\mathbb Z}

\newcommand\res{\arrowvert_}

\renewcommand{\S}{\mathcal S}
\renewcommand{\SS}{\mathcal{SS}}
\newcommand{\T}{\mathcal T}

\newcommand{\Hom}{\operatorname{Hom}}
\newcommand{\Ext}{\operatorname{Ext}}
\newcommand{\Aut}{\operatorname{Aut}}
\newcommand{\Tw}{\operatorname{Tw}}

\makeatletter \@addtoreset{equation}{section} \makeatother

\newtheorem{defn}[equation]{Definition}
\newtheorem{thm}[equation]{Theorem}
\newtheorem{lem}[equation]{Lemma}

\newtheorem{prop}[equation]{Proposition}
\newenvironment{pf}{\noindent\emph{Proof.}}{\hfill$\square$\\}
\newenvironment{rmk}{\noindent\textbf{Remark}.}{\\}

\title{\textbf{Stability conditions and the braid group}}
\author{R. P. Thomas}
\date{}

\begin{document}
\maketitle

\begin{abstract} \noindent
We find stability conditions (\cite{D2}, \cite{Br1}) on
some derived categories of differential graded modules over a graded
algebra studied in \cite{RZ}, \cite{KS}. This category arises in both derived
Fukaya categories and derived categories of coherent sheaves. This
gives the first examples of stability conditions on the A-model side of
mirror symmetry, where
the triangulated category is not naturally the derived category of an abelian
category. The existence of stability conditions, however, gives many such
abelian categories, as predicted by mirror symmetry.

In our examples in 2 dimensions we completely
describe a connected component of the space of stability conditions. It is
the universal cover of the configuration space $C_{k+1}^0$ of $k+1$
points in $\C$ with centre of mass zero, with deck transformations the braid
group action of \cite{KS}, \cite{ST}. This gives a geometric origin for these
braid group
actions and their faithfulness, and axiomatises the proposal for stability
of Lagrangians in \cite{Th} and the example proved by mean curvature flow
in \cite{TY}.
\end{abstract}


\section{Introduction}

This paper presents a result in pure algebra, but one which is motivated
entirely by geometry and physics, especially mirror symmetry. It gives
examples of Douglas and Bridgeland's notion of stability conditions for
triangulated categories \cite{D1}, \cite{D2}, \cite{AD}, \cite{Br1},
and draws together and axiomatises many of the known tests of Kontsevich's
homological mirror conjecture \cite{K} (for instance on stability of Lagrangians
in Fukaya categories \cite{Th}, Shapere-Vafa's examples of special Lagrangians
\cite{SV}, \cite{TY}, and braid groups of monodromies \cite{KS}, \cite{ST}). We
explain some of the geometry and physics background, distil this into a purely
algebraic setup, and then apply the axioms of stability conditions \cite{Br1}
to this problem. The result, which can be read independently of the previous
sections, is a description of (a connected component
of) the space of stability conditions on a certain natural triangulated
category arising in many areas of geometry and algebra \cite{RZ}, \cite{KS},
giving both an axiomatic justification
for the conjectures and results of \cite{Th}, \cite{TY}, and a geometric
``explanation" or origin for the faithful braid group actions of \cite{KS},
\cite{ST}, at least in dimension 2. There are other braid group actions in
\cite{Sz}; it would be nice to see them arise from stability conditions also.

\vspace{0.2cm}

\noindent \textbf{Acknowledgements.} This project started out as joint
work with Tom Bridgeland, who pulled out due to laziness and my misuse
of the term ``Jordan-H\"older filtration". I owe him a great deal for discussions,
an advanced copy of his axioms \cite{Br1}, and the many errors that he spotted. He has now produced a much more professional paper \cite{Br2} extending the results of this paper to all Dynkin diagrams. I would also like to thank
Ludmil Katzarkov and Paul Seidel for extremely useful conversations.
The author is supported by a Royal Society university research fellowship.


\section{Geometry and physics background}

Consider the following much-studied (\cite{KS}, \cite{SV}, \cite{TY})
affine algebraic variety $X=X^N(p)$:
\begin{equation} \label{Ak}
\left\{\sum_{i=1}^Nx_i^2=p(t)\right\} \,\subseteq\ \C^N\times\C,
\end{equation}
where $p$ is some degree $k+1$ polynomial in $t\in\C$ with only
simple zeros. There is a natural K\"ahler form $\omega$ restricted from
$\C^{N+1}$, and a nowhere-zero holomorphic volume form $\Omega_p$ (``almost Calabi-Yau" structure) given by
taking the Poincar\'e residue (\cite{GH} p 147) of the meromorphic form
$dx_1\ldots dx_{N\,}dt\big/\big(\sum x_i^2-p(t)\big)$ on $\C^{N+1}$; this
can be written as
\begin{equation} \label{form}
\Omega=\Omega_p=(-1)^{N+i+1}{dx_1\ldots\widehat{dx_i}\ldots dx_N\,dt\res{X^N}
\over2x_i}={dx_1\ldots dx_N\res{X^N}\over\dot p(t)}
\end{equation}
for any $i$ (so where $x_i=0\ \forall i$ we can use the second expression).
Here $\widehat{dx_i}$ means that we omit the $dx_i$ term from the wedge product.

Each smooth fibre over $t\in\C$ is an affine quadric with a natural
Lagrangian $S^{N-1}$
``real" slice, namely the intersection of the fibre with the slice
$$
x_i\in\sqrt{p(t)}\,\R\quad\forall i.
$$
It is invariant under the obvious $O(N)$ action on $X^N$, and is the
vanishing cycle of every singular fibre (i.e.\ the fibres over the
roots of $p$). Therefore any path
$\gamma:\,[0,1]\to\C$ from one zero of $p$ to another (avoiding zeros of $p$
in its interior) lifts to give a canonical
$O(N)$-invariant Lagrangian $N$-sphere, $S^{N-1}$-fibred
over $\gamma$ except at the endpoints where it closes up.

Picking an ordering of the $k+1$ zeros of $p$, we can pick $k$ consecutive
paths in $\C$ joining one zero to the next, forming a so-called ``$A_k$-chain" of paths. The corresponding Lagrangian $N$-spheres $\{L_i\}_{i=1}^k$ in $X^N$
(the vanishing cycles of the $A_k$-singularity; the singularity
obtained by putting $p(t)=t^{k+1}$ in (\ref{Ak})) form a basis for its
homology and define an $A_k$-chain of Lagrangians \cite{Se1}. That is, their
\emph{geometric} intersections satisfy
$$
|L_i.L_j|=1,\ \ |i-j|=1, \quad\qquad |L_i.L_j|=0,\ \ |i-j|>1.
$$
Choosing their gradings \cite{Se2} appropriately, we can arrange
that their Floer cohomologies satisfy
\begin{eqnarray} \nonumber
HF^0(L_i,L_i) \!&=&\! \C\ =\ HF^N(L_i,L_i), \\
HF^1(L_i,L_{i+1}) \!&=&\! \C\ =\ HF^{N-1}(L_{i+1},L_i), \label{sph}
\end{eqnarray}
for all $i$, with all other groups zero.

As $p$ varies, all of these varieties $X_p$ are \emph{isomorphic as symplectic
manifolds}, even though as complex manifolds they are varying. Thus around
loops in the space of polynomials with simple zeros $\{p\}$ we get monodromy in the symplectic automorphism group $\Aut(X,\omega)$. Scaling $p$ makes
no difference to the symplectic geometry, nor, in dimension $N=2$ to the
\emph{graded} symplectic geometry, as $\Omega_p$ (\ref{form}) is left unaltered.
Similarly everything is invariant under translations in the $t$-plane; dividing
by these reparametrisations leaves us with the (simple) roots
of $p$ normalised (say) to have centre of mass zero. That is we get
the configuration space $C_{k+1}^0$ of $k+1$ points of mean zero in $\C$,
with fundamental group
the braid group $B_k$. This gives a representation of $B_k\to\Aut(X,\omega)$,
generated by ``generalised Dehn twists" about the vanishing cycles; this
amounts to ``twisting" about the Lagrangian fibred over a path between two
roots of $p$ by rotating the two points in $\C$ about each other; i.e.\ it
arises from the usual description of $B_k$ lifted from $\C$ (with $k+1$ marked
points) to $X$. In fact these Dehn twists give such a braid group action
on any symplectic manifold with an $A_k$-chain of Lagrangian spheres \cite{Se1},
and so this a local model for symplectic automorphisms of compact symplectic
manifolds, in particular Calabi-Yau manifolds. This induces a braid group
action on the derived Fukaya category $D^b(Fuk(X))$ of $X$ (though one
must check first, as in \cite{Se2}, that the symplectomorphisms lift naturally
to the graded symplectomorphism group).

Under mirror symmetry (and there is a proposal for the mirror of this
$A_k$-smoothing in \cite{ST} Section 3f)
we cannot expect a braid group action of holomorphic
automorphisms of any mirror complex manifold $Y$; there are in general very
few holomorphic automorphisms of Calabi-Yau manifolds, and it is one of
Kontsevich's great insights that the mirror is really the bounded derived
category $D^b(Y)\cong D^b(Fuk(X))$ of coherent sheaves
on the Calabi-Yau rather than the Calabi-Yau itself; in this way the
automorphisms can be matched, so things work in families.

\begin{figure}[h]
\center{
\input{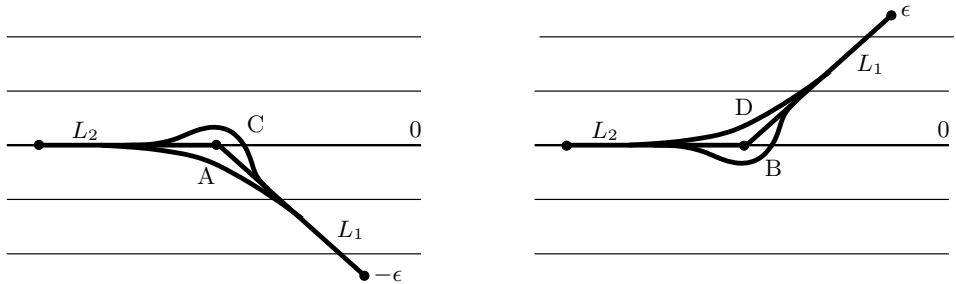}
\caption{The connect sums $L_1\#L_2$ (A -- stable, B -- unstable) and
$L_2\#L_1$ (C -- unstable, D -- stable) in 2-dimensions, in 2 different
complex structures. \label{stab2}}}
\end{figure}

So if the only vestiges of our geometric picture above that remain under
the mirror map are categorical, we would like to see the varying complex
structure (given by the polynomial $p$) that induces the braid monodromy
at the level of the (purely symplectically defined) derived Fukaya category.
How one ``ought" to do this has long be known to physicists, and was described
to me many years ago by Eric Zaslow. Namely one should follow
the D-branes in $D^b(Fuk(Y))$; these depend on the complex structure, and
are thought to be the \emph{special Lagrangians} (Lagrangians $L$ for which
the $N$-form (\ref{form}) $\Omega|_L$ has constant phase on $L$).
On going round a loop in complex structure
moduli space the set of special Lagrangians undergoes monodromy which one might
hope extends uniquely to the full derived category. Similarly in the
mirror picture the D-branes in $D^b(X)$ depend on the K\"ahler structure,
and are to a first approximation the \emph{stable} bundles that satisfy
an appropriately perturbed Hermitian-Yang-Mills equation. In \cite{Th} this
picture was mirrored to give an appropriate notion of stability for Lagrangians
which should conjecturally be equivalent to the Lagrangian being a hamiltonian
deformation of a (unique) special Lagrangian. In the case of the Lagrangians of
our above example, this can all be illustrated very simply \cite{TY}, at least in
$N=2$ dimensions (higher dimensions are more complicated, but also illustrated
in \cite{TY}). $O(N)$-invariant \emph{graded} \cite{Se2}
Lagrangian spheres correspond to paths in $\C$ between distinct zeros of
$p(t)$, missing other zeros, with a continuous choice of lift to
$\R$ of the phase $\in\R/2\pi\Z$ of the tangent to the path at each point.
Special Lagrangians correspond to (spheres fibred over)
straight lines, and $O(N)$-invariant hamiltonian isotopies to isotopies
of the path in $\C$, with endpoints fixed, not crossing
any other zeros of $p$. Extensions correspond to graded Lagrangian connect
sums $\#$, and stability of $L$ to there being no such connect sum $L_1\#L_2$
hamiltonian isotopic to $L$, with the phase of $L_1$ greater than that
of $L_2$; see Figure \ref{stab2}, in which $L_1$ has phase $\pm\epsilon$
and $L_2$ has phase zero.

Simple examples of loops of complex structures $p$ do indeed exhibit
the braid monodromy action on the set of such special Lagrangians \cite{TY}.
For a long time it was not clear how to extend such naive notions of stability
(in terms of injective morphisms or subobjects) to full triangulated categories
(where there is no notion of subobject). But the beautiful work of Douglas
\cite{D2}, now axiomatised in \cite{Br1}, purports to give the answer; the
purpose of this paper is to show how it works in our $A_k$ example, and how
it agrees with the naive definitions mentioned above. This also gives a
new example where the axioms of \cite{Br1} hold; in general examples are hard
to find. Known examples work on the B-model side where there is an abelian
category (the coherent sheaves), and where the usual notions
of (semi)stability of sheaves are well understood and give examples satisfying
the axioms (with the semistable objects of the axioms below being the semistable
sheaves and their translations) which can be deformed to give more.
In the example of this paper there is no natural abelian
category, as we start on the A-model side with the Fukaya category. Using
the stable objects of \cite{TY} we prove the existence of stability conditions,
which in turn shows that the derived Fukaya category is the
derived category of many abelian categories.

Our results also apply to the B-model derived category of coherent sheaves,
and thus gives a geometric, monodromy background for the mirror braid group action
on the derived category described in \cite{ST}. The varying complex structure
$p$ is then interpreted as a varying complexified K\"ahler form on the
mirror manifold containing an $A_k$-chain of spherical objects. In fact
we do not use $p$, but the closely related set of the values it gives to
the $A_k$-chain of Lagrangians via (\ref{form}):
\begin{equation} \label{addlags}
\left\{\sum_{j=1}^iZ(L_j)=\int_{L_1\cup\ldots\cup L_i}\!\!
\Omega_p\right\}_{i=1}^k.
\end{equation}
In two dimensions $Z(L_i)$ is a constant times by the vector in $\C$
between the two zeros of $p$ that form the endpoints of the path over which
$L_i$ is $S^1$-fibred (see (A.1) of \cite{Ho}). Thus, up to an additive constant
in $\C$ which does not affect the homotopy type, the numbers (\ref{addlags}) and the origin together constitute the distinct roots of $p$ in $\C$, and the map $p\mapsto\big\{0,\sum_{j\le i}Z(L_j)\big\}$ is a $\C^*$-bundle inducing an isomorphism $\pi_1/\pi_1(\C^*)\to\pi_1(C_{k+1})\cong B_k$.
(In higher dimensions the relationship
is more complicated; for small paths, $Z(L_i)$ is roughly the
$(n/2)$th power of the vector in $\C$ representing $L_i$, which is why
it is simpler to stick to two dimensions for our analysis of monodromy.)
Instead of the space of $p$\,s, then, we will find the space of stability
conditions to be the universal cover of the configuration space $C_{k+1}^0$
of $(k+1)$ distinct points of mean zero in $\C$, with fundamental group $B_k$.

Most of the work can be done
by just dealing with curves in the plane (with endpoints in the $k+1$ marked
points) as in \cite{TY}, and the stability conditions that emerge axiomatise
that of \cite{Th},
\cite{TY}; i.e.\ relate to whether or not the path can be pulled ``tight"
(straight, in the dimension $N=2$ case, or to a special Lagrangian of \cite{SV}
in general; this was done in \cite{TY} by mean curvature flow) without
touching one of the other marked points. Unfortunately one has to
deal with slightly more general objects in $D^b(Fuk(X))$
than can be represented solely by curves in the plane. This is because
the curves do not
form an abelian, or triangulated, category -- in general one cannot form
the kernel or cokernel or cone on a morphism (element of Floer homology)
between two Lagrangians. They do form an $A^\infty$ category in a complicated
way (which involves difficult counting of holomorphic discs with boundaries
in the Lagrangians), though an intrinsic formality result \cite{ST} means
we will not have to worry about the precise $A^\infty$ structure.
We then have to derive this category, a formal procedure producing cones on
abstract morphisms that introduces extra objects not all representable by curves.
As proposed in \cite{Th}, $\Ext^1$s are geometrically represented by (graded) connect
sums of Lagrangians; what stops us from using this to form geometric representatives
of all cones is the fact that some of the $\Ext^1$s lead to immersed Lagrangians,
whose Floer homology is not well understood, and some lead to different representatives
of the same class in the Fukaya category if one either takes the relative
connect sum \cite{Th} or perturbs and takes a transverse connect sum. It
would be nice to find a purely one-dimensional geometric description of
all objects of this category.

Nonetheless our constructions are motivated by pictures for those objects
of $D^b(Fuk)$ that \emph{are} representable by curves, and so such figures
often accompany algebra below without any explanation; these should
be helpful to anyone who has read the above background. Unfortunately such
physical arguments do not suffice for the whole category for the $A^\infty$
reasons mentioned above, but nonetheless give a good intuitive
picture for the axioms \cite{Br1} that will be familiar to physicists, reminiscent
of ``marginal stability".

To show that our results are relevant to ``real life" we now define our
category and show it is indeed faithfully included in the derived
Fukaya category of any symplectic manifold with an $A_k$-chain of Lagrangian
spheres, and in the derived category of any smooth projective variety with
an $A_k$-chain of spherical objects \cite{ST}.

\section{The categories}

We start by defining a simple auxiliary category to define our graded algebra
(closely related to the one studied in \cite{RZ}, \cite{KS}) whose derived category
of dgms will be the triangulated category of study. Compare (\ref{sph}).

\begin{defn} \label{cat}
Let $\mathcal C=\mathcal C^N_k$ be the $\C$-linear graded category 
of an $A_k$-chain $\{E_j\}_{j=1}^k$ in dimension $N$; the category with $k$
\emph{spherical} objects $E_j$, \vspace{1mm} \\
$\bullet\ \Hom^i(E_j,E_j)=\left\{\begin{array}{cl} \C & i=0,N \\ 0 &
\mathrm{otherwise,} \end{array} \right.$ \vspace{1mm} \\
and morphisms \vspace{1mm} \\
$\bullet\ \Hom^i(E_j,E_{j+1})=\left\{\begin{array}{cl} \C & i=1 \\ 0 &
\mathrm{otherwise} \end{array} \right.$ \\
$\bullet\ \Hom^i(E_{j+1},E_j)=\left\{\begin{array}{cl} \C & i=N-1 \\ 0 &
\mathrm{otherwise} \end{array} \right.$ \\
$\bullet\ \Hom^i(E_j,E_k)=0,\ \ |j-k|>1.$ \vspace{2mm} \\
Denoting by $(j,j\!\pm\!1)$ the generator of $\Hom^i(E_j,E_{j\pm1})$,
we impose the relation that $(j,j\!\pm\!1)\circ(j\!\pm\!1,j)$ is the generator
$f_j$ of $\Hom^N(E_j,E_j)$ for all $j$; in particular then $(j,j+1)(j+1,j)=
(j,j-1)(j-1,j)=f_j$, and the category has a duality: the pairing $\Hom^i(F,E)
\otimes\Hom^{N-i}(E,F)\to\Hom^N(E,E)\cong\C$
is perfect for all $E,F\in\mathcal C$.
\end{defn}

\begin{figure}[hb]
\begin{center}
\epsfig{file=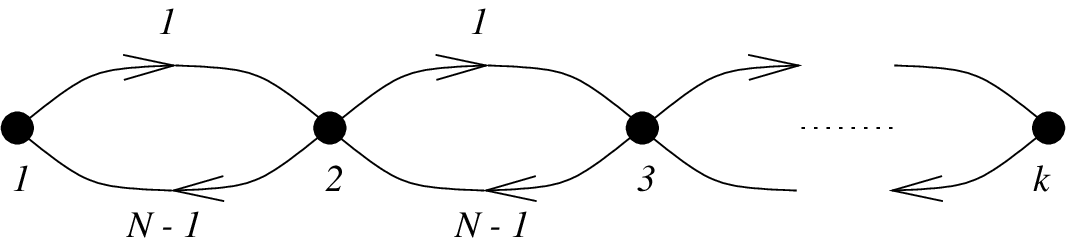} \\ \caption{\label{quiver}}
\end{center}
\end{figure}

In the usual way $\mathcal C$ defines a graded unital algebra End\,$\mathcal C$:
$$
A_k=A^N_k=\mathrm{End\,}\big(\bigoplus_{j}E_j\big)=
\bigoplus_{ijk}\Hom^i(E_j,E_k).
$$
This is a quotient of the path algebra of the quiver in Figure \ref{quiver}.
Generators, over $\C$, of the algebra are given by oriented paths, graded
as in the Figure, with multiplication given by composition of paths (this
is the \emph{dual} picture to the Lagrangian one described in the introduction,
with vertices in Figure \ref{quiver} corresponding to Lagrangians, i.e.\
paths between zeros of $p$, and arrows corresponding to intersections of
Lagrangians, i.e.\ zeros of $p$). The quotient is by
the two sided ideal generated by $(i,i-1,i)-(i,i+1,i)$, $(i-1,i,i+1)$ and
$(i+1,i,i-1)$ for all $i=2,\dots,k-1$, in the obvious notation; see (\cite{ST}
Section 4c) for more details. Denote by $e_i$ the identity in $\Hom(E_i,E_i)$
(i.e.\ the constant path at the $i$th node), so that $1=\sum_ie_i$ and left
(right) multiplication by $e_i$ is projection onto those paths that begin
(end) at the $i$th node. Similarly we denote $f_i=(i,i-1,i)=(i,i+1,i)$
(or whichever of the two is defined, if $i=1$ or $k$).

Now form the bounded derived category $D(A_k)$ of differential graded (right)
modules (dgms) over $A_k$; again see (\cite{ST} Section 4a).
This is \emph{not} the derived category of the abelian category of dgms
over $A_k$, but the localisation of this category by quasi-isomorphisms
(dgm maps that induce isomorphisms on the cohomology graded modules). It
is, however, triangulated. Denote by $P_i$ the projective (right) module
$e_iA_k$; these form an $A_k$-chain in $D(A_k)$ in the sense of
Definition \ref{cat} (with $\Hom^i$ the $i$th cohomology of Hom)
and so form a subcategory whose cohomology category is isomorphic to the
original $\mathcal C_k$.

\begin{defn} \label{Dk}
$D_k=D_k^N$ is the sub-triangulated category of $D(A_k)$ generated by the
$P_i$; the smallest triangulated subcategory containing the $P_i$.
\end{defn}

Note that the duality (\ref{cat}) induces a ``trace map" $\Hom^N(E,E)\to\C$
for all $E\in D_k$ such that there exists a perfect pairing $\Hom^i(F,E)
\otimes\Hom^{N-i}(E,F)\to\Hom^N(E,E)\to\C$, i.e.\ a duality
\begin{equation} \label{dual}
\Hom^{N-i}(E,F)^*\cong\Hom^i(F,E),
\end{equation}
functorial in $E$ and $F$, for all $E,F\in D_k$.

This is the category whose stability conditions we will study in this paper.
Before we do, we point out that due to the intrinsic formality result of
\cite{ST}, it is contained in the derived
Fukaya category of a symplectic manifold containing a (suitably graded)
$A_k$-chain of Lagrangian spheres, and, mirror to this, in any derived
category of coherent sheaves containing a (suitably shifted) $A_k$-chain
of spherical objects $\{E_i\}$ \cite{ST}. Many thanks to Paul Seidel for
this argument.

In the first case, the full
subcategory $\mathcal A$ defined by the $A_k$-chain defines a natural
dga whose cohomology algebra is our $A^N_k$. By intrinsic formality \cite{ST},
then, $\mathcal A$ is actually quasi-isomorphic to $A^N_k$.
Since the derived Fukaya category contains cones, $\Tw(\mathcal A)$,
the triangulated category of twisted complexes \cite{BK} on $\mathcal A$,
is a full subcategory of the derived Fukaya
category. But it is quasi-equivalent to $\Tw(A^N_k)$, which is equivalent
to the derived category of dgms generated by the projective modules $P_i=e_iA$
above (since the functor $\Hom(\ \cdot\ ,\oplus_iE_i)$ from $C^N_k$ to
the category of dgms over $A^N_k=\Hom(\oplus_iE_i,\oplus_iE_i)$ is full
and faithful by the Yoneda lemma, and takes the objects $E_i$ of $C^N_k$ to
the projectives $P_i$).

Similarly, given an $A_k$-chain in a derived category of coherent sheaves
over a smooth projective variety, we can work with the equivalent homotopy
category of complexes of quasicoherent injective sheaves
with bounded coherent cohomology. This is a dg-category containing cones,
so is isomorphic to its own Tw. The $A_k$-chain defines a full sub-dg-category
defining a dga quasi-isomorphic to its cohomology graded algebra $A^N_k$
by intrinsic formality. The rest of the argument is then the same.

\begin{prop} For $N\ge2$, $D^N_k$ is fully faithfully
embedded in the derived Fukaya category
of any $2N$-dimensional symplectic manifold containing an $A_k$-chain of spheres;
similarly for the derived category of coherent sheaves on any smooth quasiprojective
$N$-dimensional variety with an $A_k$-chain of spherical objects. \hfill$\square$
\end{prop}

For instance in dimension two we can consider the coherent
sheaves over $\C^{\,2}$ (finitely generated modules over $\C\,[x,y]$) supported at
the origin, and the standard $SU(2)$ action of $\mathbb Z/k$ on
$\C^{\,2}$. Then the derived category $D^{\Z/k}_0(\C^{\,2})$ of equivariant sheaves
supported at the origin is equivalent to the derived category of coherent
sheaves on the minimal resolution $\widehat{\C^{\,2}/(\Z/k)}$ supported on the
exceptional set $E$ \cite{KV}, \cite{BKR}. The exceptional locus $E$
is an $A_k$-chain of $-2$-spheres whose structure sheaves form an $A_k$-chain
in the derived category \cite{ST}, and correspond under the equivalence
to the non-trivial irreducible representations of $\Z/k$. So $D^2_k$ is
embedded in $D^{\Z/k}_{0}(\C^{\,2})\cong D^b_E\left(\widehat{\C^{\,2}
\over\Z/k}\right)$.

\section{Stability}

In \cite{Br1}, a notion of stability for derived categories is given, axiomatising
the proposal of Douglas \cite{D2}. Let $K(\T)$ denote the Grothendieck
group of $\T$.

\begin{defn} \label{def}
A \emph{stability condition} $(Z,\SS)$ on
a triangulated category $\T$ consists of a
linear map $Z\colon K(\T)\to\C$
and full subcategories
$\SS(\phi)\subset\T$ for each $\phi\in\R$ satisfying the following five axioms:

\begin{itemize}
\item[(a)] for all $\phi\in\R$, $\SS(\phi+1)=\SS(\phi)[1]$,

\item[(b)] if $E\in \SS(\phi)$ then $Z(E)=m(E)\exp(i\pi\phi)$ with $m(E)>0$,

\item[(c)] for $0\neq E\in\T$ there is a finite sequence of real
numbers
\[\phi_1>\phi_2> \cdots >\phi_n\]

and a collection of triangles
\[
\xymatrix@C=.5em{
0_{\ } \ar@{=}[r] & E_0 \ar[rrrr] &&&& E_1 \ar[rrrr] \ar[dll] &&&& E_2
\ar[rr] \ar[dll] && \ldots \ar[rr] && E_{n-1}
\ar[rrrr] &&&& E_n \ar[dll] \ar@{=}[r] &  E_{\ } \\
&&& S_1 \ar@{-->}[ull] &&&& S_2 \ar@{-->}[ull] &&&&&&&& S_n \ar@{-->}[ull] 
}
\]
with $S_i\in\SS(\phi_i)$ for all $i$,

\item[(d)] if $\phi_1>\phi_2$ and $S_i\in\SS(\phi_i)$ then $\Hom_{\T}(S_1,S_2)=0$,
\item[(e)] the subset
\[Z\,\Big(\bigcup_{\phi\in\R} \SS(\phi)\Big)\subset\C\]
has no limit points in $\C$.
\end{itemize}
\end{defn}

The map $Z$ is called the \emph{central charge} of the stability condition.
The objects of the subcategory $\SS(\phi)$ are said to be \emph{semistable
of phase $\phi$}; the simple semistables are \emph{stable}, and we denote
these by $\S(\phi)$. We call the choice
of a lift of the phase $e^{i\pi\phi}$ of such an object $E$ to a real number
$\phi$ a \emph{grading} of $E$. It is an easy exercise to
check that the decomposition of a nonzero object $E$ given in (c) is unique;
the objects $S_i$ are called the \emph{semistable factors} of $E$. We sometimes
call such a collection of triangles a filtration. The \emph{mass} of $E$ is
the positive real number
\begin{equation} \label{mass}
m(E)=\sum_i|Z(S_i)|.
\end{equation}
By the triangle inequality one has $m(E)\geq|Z(E)|$
with equality if $E$ is semistable.

The axioms are modelled on semistability for sheaves on complex curves:
on filtering objects by their cohomology sheaves, and these in turn by
their Harder-Narasimhan filtrations, we get an example satisfying axiom (c)
above.

To make this definition more manageable we give some conditions that will
imply the difficult axiom (c).

\begin{thm} \label{split}
Given $(\T,Z,\S(\phi))$, and defining $\SS(\phi)$ to consist
of all possible extensions of elements of $\S(\phi)$, suppose that these
satisfy axioms (a,b,d,e)
above. Suppose that $\bigcup_\phi\S(\phi)=\{S_i[m]:i=1,\cdots,k,\ m\in\Z\}$,
for some finite set of $S_i$ which generate $\T$
(i.e.\ every object in $\T$ is a finite extension of shifts of $S_i$s). Suppose
also that for any non-trivial element of $\Ext^1(F,E)$ defining a triangle
\begin{equation} \label{ext}
E\to C\to F, \quad \text{with}\ \ E\in\S(\phi),\ F\in\S(\psi),\ \ \phi<\psi,
\end{equation}
$C$ is either in $\S(\theta)$
with $\theta\in(\phi,\psi)$, or $C=A\oplus B$, with
$A\in\S(\alpha),\,B\in\S(\beta),\ \phi<\alpha\le\beta<\psi$.
Then $(\T,Z,\SS(\phi))$ satisfy axiom (c), i.e.\ they define a stability
condition.
\end{thm}

\begin{rmk}
This theorem is surely true more generally: that we get a stability condition
if all extensions (\ref{ext}) have a Harder-Narasimhan filtration by shifts
of $S_i$s (rather than just splitting into two stable objects),
but we will only require the above result.
\end{rmk}

\begin{pf}
We have to find a Harder-Narasimhan filtration (c) for any object $E\in\T$.
Since the $S_i$ generate $\T$, we can find a collection of triangles
\[
\xymatrix@C=.5em{
0_{\ } \ar@{=}[r] & E_0 \ar[rrrr] &&&& E_1 \ar[rrrr] \ar[dll] &&&& E_2
\ar[rr] \ar[dll] && \ldots \ar[rr] && E_{n-1}
\ar[rrrr] &&&& E_n \ar[dll] \ar@{=}[r] &  E_{\ } \\
&&& Q_1 \ar@{-->}[ull] &&&& Q_2 \ar@{-->}[ull] &&&&&&&& Q_n \ar@{-->}[ull] 
}\] \vspace{-11mm}
\begin{equation} \label{HN} \vspace{2mm}
\end{equation}
with each $Q_j\in\bigcup_{i,m}\{S_i[m]\}$ \emph{stable} of phase $\phi_j$.

Suppose that for some $i$, $\phi_i<\phi_{i+1}$. Replace $E_{i-1}\to E_i\to E_{i+1}$
in the above filtration by the composition $E_{i-1}\to E_{i+1}$; this has
cone $Q$ fitting into a triangle
\begin{equation} \label{Q}
Q_i\to Q\to Q_{i+1}.
\end{equation}
We can now use the assumption on such extensions, as $\phi_i<\phi_{i+1}$
and the $Q_i$ are stable.

Either (i) $Q$ is stable, and we replace (\ref{HN}) by our new filtration with
one less triangle (with $Q_i,\,Q_{i+1}$ replaced by $Q$ of phase
$\phi\in(\phi_i,\phi_{i+1})$). We then start the process again,
looking for $j$ such that $\phi_j<\phi_{j+1}$.

Or (ii) $Q=0$, so $Q_{i+1}=Q_i[1]$; in this case we also remove $E_{i+1}$
from the filtration to give a new filtration with $E_{i-1}\to E_{i+2}$
in the middle (with cone $Q_{i+2}$ forming the new triangle).

Or (iii) $Q=Q_i'\oplus Q_{i+1}'$ is the direct sum of two stables
of phases $\phi_i'\ge\phi_{i+1}'$ (without loss of generality) in
$[\phi_i,\phi_{i+1}]$ (with the closed interval being necessary in case
the extension (\ref{Q}) is the trivial one). In this case we define
$E'_i$ by the triangle
$$
E'_i\to E_{i+1}\to Q'_{i+1},
$$
where the second arrow is the composition of $E_{i+1}\to Q$ and $Q=Q_i'\oplus Q_{i+1}'
\to Q_{i+1}'$. Then a small check with the octahedral Lemma gives
us a new filtration (\ref{HN}) with $E_{i-1}\to E_i\to E_{i+1}$ replaced
by
\[
\xymatrix@C=.5em{
E_{i-1} \ar[rrrr] &&&& E'_i \ar[rrrr] \ar[dll] &&&& E_{i+1}
\ar[dll] \\
&& Q'_i \ar@{-->}[ull] &&&& Q'_{i+1} \ar@{-->}[ull]
}\]
Again we now start again with this filtration.

We claim this procedure terminates; that is after a finite number of steps
we have that $\phi_i\ge\phi_{i+1}$ for all $i$. To demonstrate this, we
assign to each such filtration a number which is both bounded below and
decreases (by some bounded below amount) at each stage. Firstly, we may
assume without loss of generality (by replacing $E$ by $E[k]$ for some
$k$ if necessary) that each phase $\phi_i$ of the $Q_i$ in (\ref{HN}) is
positive.

Then to (\ref{HN}) we can associate the real number $\sum_{k=1}^nf(k)\phi_k$
for some strictly positive, strictly increasing, concave function $f(x)$.
Then in case (i) above this clearly decreases, as $\phi<\phi_{i+1}$,
so that $f(i)\phi<f(i+1)\phi_{i+1}<f(i)\phi_i+f(i+1)\phi_{i+1}$, and the
sum over all higher $k\ge i+1$ is also smaller.
Case (ii) is even clearer (remembering that all $\phi_i>0$). Finally for
case (iii) we pick a sufficiently concave function $f$ such that
$$
f(x)\phi+f(x+1)\psi>f(x)\beta+f(x+1)\alpha,
$$
for all $x\ge0$ and $\phi<\alpha\le\beta<\psi$ coming from phases of extensions of stable objects as in the assumptions of the Theorem; this is possible since the number of such extensions is finite (up to shifts, which leave the above inequality unaffected).

Then in case (iii) $f(i)\phi_i+f(i+1)\phi_{i+1}>f(i)\phi_i'+f(i+1)\phi_{i+1}'$
(i.e.\ the above inequality with $\phi_i'=\beta,\,\phi_{i+1}'=\alpha$) ensures
that the functional again decreases.

This procedure now terminates as the amount the functional decreases by is bounded below by the discreteness of the phases of the stable objects.
This gives us a ``Jordan-H\"older filtration" of $E$ into stable objects
of nonincreasing phases (we use the term as in holomorphic bundle theory
or \cite{TY} for Lagrangians, rather than the category theoretic terminology).
To get the Harder-Narasimhan filtration (c) we
bundle together any $Q_i$s of the same phase. That is, if $\phi_i=\phi_{i+1}$,
replace the $E_{i-1}\to E_i\to E_{i+1}$ part of the filtration by just
$E_{i-1}\to E_{i+1}$, with one less triangle with cone $Q$ fitting into
a triangle
$$
Q_i\to Q\to Q_{i+1}.
$$
By assumption this means $Q$ is semistable of phase $\phi_i=\phi_{i+1}$;
now continue the process until the phases are strictly decreasing.
\end{pf}

We mention in passing that given a stability condition we get
a family of bounded nondegenerate t-structures $\mathcal F_t=\mathcal
F_{t+1}[-1]$ on $\T$ given
by the full subcategory of objects whose semistable factors (c) all of
phase $>t$ \cite{Br1}. This has as its heart
\begin{equation} \label{heart}
\mathcal F_t\cap \mathcal F_{t+1}^\perp,
\end{equation}
the full subcategory of objects whose semistable
factors (c) all have phase in $(t,t+1]$, so in particular we can assign
a phase $\phi(A)\in(t,t+1]$ to each object $A\ne0$ of the heart such
that $Z(A)=m(A)e^{i\pi\phi(A)}$ with $m(A)>0$. It turns out then \cite{Br1}
that $A$ is semistable if and only if for every subobject $B$ of $A$ in
the heart (which is, recall, an abelian category), $\phi(B)\le\phi(A)$. We
may also take $B$ to be stable in this test. We
also get (a more standard) Harder-Narasimhan filtration (c) of objects
in the heart by semistable objects of the heart.

\section{Our example}

We first need some technical results about our projective modules $P_i$ in $D^N_k$.
Given modules $A,\,B$ and an element $e$ of $\Ext^1(B,A)=\Hom^0(B[-1],A)$, we will
often denote by $A\#B$ the corresponding extension (i.e.\ the cone on $B[-1]\to
A$) fitting into the triangle
$$
B[-1]\to A\to A\#B\to B.
$$
This defines a canonical $\tilde e\in\Hom^0(A\#B,B)$.

\begin{lem} \label{con}
Given $e\in\Ext^1(C,A)$ and $f\in\Ext^1(C,B)$ defining $A\#C,\,B\#C$ and
$\tilde e\in\Hom^0(A\#C,C),\,\tilde f\in\Hom^0(B\#C,C)$, form
$e\cup\tilde f\in\Ext^1(B\#C,A)$ and $f\cup\tilde e\in\Ext^1(A\#C,B)$.
Then the corresponding extensions are isomorphic:
$$
A\#(B\#C)\cong B\#(A\#C).
$$
Similarly we have
$$
A\#(B\#C)\cong (A\#B)\#C,
$$
if the classes in $\Ext^1(B\#C,A),\ \Ext^1(C,A\#B)$ defining them map to
the classes in $\Ext^1(B,A),\ \Ext^1(C,B)$ defining $A\#B,\,B\#C$ respectively.
\end{lem}

\begin{pf}
The first assertion follows from the diagram of triangles
$$
\xymatrix@R=.9em{
& B \ar[d] \ar@{=}[r] & B \ar[d] \\
A \ar[r] \ar@{=}[d] & A\#(B\#C) \ar[r] \ar[d] & B\#C \ar[d] \\
A \ar[r] & A\#C \ar[r] & C\,.\!\!
}$$
Here we start with the bottom two rows and left and right hand columns,
and this defines the arrow $A\#(B\#C)\to A\#C$ by taking cones. The octahedral
Lemma then gives the top row, so the central column now shows that $A\#(B\#C)\cong
B\#(A\#C)$.

The second statement follows from similar yoga around the diagram
$$ \vspace{-6mm}
\xymatrix@R=.9em{
A \ar[r] \ar@{=}[d] & A\#B \ar[r] \ar[d] & B \ar[d] \\
A \ar[r] & E \ar[r] \ar[d] & B\#C \ar[d] \\
& C \ar@{=}[r] & C
}$$
\end{pf}

\begin{figure}[h]
\center{
\input{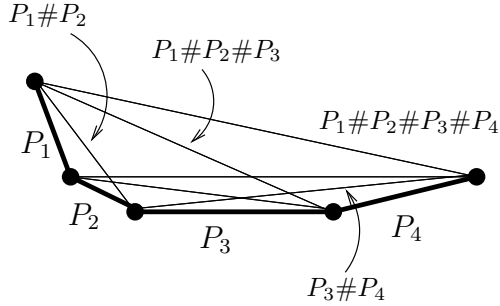}
\caption{Our collection of stable objects \label{slag}}}
\end{figure}

In our use of this Lemma below, $\#$ will be the \emph{unique} nontrivial
extension between the objects concerned. We use the notation $\cong$ for
quasi-isomorphism and $\C\,[-n]$ for a copy of $\C$ shifted into degree $n$,
so that $\Hom(E,F)\cong\C\,[-n]$ is equivalent to $\Hom^*(E,F)=\C$ for $*=n$
and zero otherwise. The following is best interpreted in terms of pictures
such as Figure \ref{slag} and the discussion of \emph{graded} connect sums
and relative connect sums in \cite{Th}, \cite{TY}.

\begin{prop} Define $P_{ii}:=P_i$, then inductively (on $j\ge i$) one can define
$P_{i,j+1}:=P_{ij}\#P_{j+1}$ by $\Hom^1(P_{j+1},P_{ij})=\C$. \\

\noindent Moreover we then have $\,\Hom(P_k,P_{ij})\cong\left\{\!\!\begin{array}{ll}
\C\,[1\!-\!N] & k=i-1 \\ \C & k=i \\ \C\,[-N] & k=j \\ \C\,[-1] & k=j+1 \end{array}
\right.$ and zero otherwise.
\end{prop}

\begin{pf}
The result is true for $P_{ii}$, i.e.\ $j=i$, by the fact that the $P_i$
form an $A_k^N$-chain of spherical objects (\ref{cat}). Inductively then,
assume it is true for $j$, so that we have defined $P_{i,j+1}=P_{ij}\#P_{j+1}$.

Then the triangle
$$
P_{ij}\to P_{i,j+1}\to P_{j+1}\to P_{ij}[1]
$$
gives
\begin{equation} \label{homs}
\Hom(P_k,P_{ij})\to\Hom(P_k,P_{i,j+1})\to\Hom(P_k,P_{j+1})\to\Hom(P_k,P_{ij})[1].
\end{equation}

\noindent $\bullet$ For $k<i-1,\,k>j+2,$ or $i<k<j$, the second complex is
(quasi-isomorphic to) zero since all the others
are by the induction assumption.

\noindent $\bullet$ For $k=i-1$ (respectively $k=i$) (\ref{homs}) becomes
$$
\C\,[-d]\to\Hom(P_k,P_{i,j+1})\to0\to\C\,[1-d],
$$
where $d=N-1$ (respectively $d=0$), so that $\Hom(P_k,P_{i,j+1})\cong\C\,[-d]$
as required.

\noindent $\bullet$ When $k=j$ (\ref{homs}) becomes
$$
\C\,[-N]\to\Hom(P_j,P_{i,j+1})\to\C\,[1-N]\to\C\,[1-N].
$$
We claim this last map $\Hom^{N-1}(P_j,P_{j+1})\to\Hom^N(P_j,P_{ij})$ is an
isomorphism,
so that $\Hom(P_j,P_{i,j+1})\cong0$ as required. This follows from the fact
that the composition to $\Hom^N(P_j,P_j)=\C$, is, by construction, the
Yoneda product
$$
\Hom^{N-1}(P_j,P_{j+1})\otimes\Hom^1(P_{j+1},P_j)\to\Hom^N(P_j,P_j),
$$
which is an isomorphism by the duality (\ref{dual}).

\noindent $\bullet$ For $k=j+1$ we get
$$
\C\,[-1]\to\Hom(P_{j+1},P_{i,j+1})\to\C\oplus\C\,[-N]\to\C,
$$
with the $\Hom^0(P_{j+1},P_{j+1})\to\Hom^1(P_{j+1},P_{ij})$ component of
the last map takes the identity to the extension class defining
$P_{i,j+1}=P_{ij}\#P_{j+1}$. Since this was chosen to be nontrivial, this
map is an isomorphism, so $\Hom(P_{j+1},P_{i,j+1})\cong\C\,[-N]$.

\noindent $\bullet$ Finally, when $k=j+2$, we have by (\ref{homs})
$$
0\to\Hom(P_{j+2},P_{i,j+1})\to\C\,[-1]\to0,
$$
as required. In particular, this now allows us to define $P_{i,j+2}$ and
continue the induction.
\end{pf}

Thus $\Hom^1(P_{ij},P_{i-1})=\C$ and we can form $P_{i-1}\#P_{ij}$, which
by inductive use of Lemma \ref{con} is $P_{i-1,j}$. (Alternatively,
$\Hom^0(P_{i-1},P_{i-1,j})=\C$ with cone $P_{ij}$ gives the same result,
as the extension cannot be trivial by the simplicity of $P_{i-1,j}$ demonstrated
in Proposition \ref{morehoms} below.) More generally, there is a unique extension
$P_{ij}\#P_{jk}\cong P_{ik}$ (again see Proposition \ref{morehoms}, for instance).
That is, we may write
$P_{ij}=P_i\#P_{i+1}\#\ldots\#P_j$ without confusion.

To define a stability condition on our category $D_k^N$, we need to define
$Z(P_i)$ for all $i$, and the set of (semi)stable objects. We assume that $N\ge2$,
so that there are no nontrivial extensions
between any $P_i$ and itself. Fix any sequence of
positive real numbers $m_i$, and a sequence of real numbers
\begin{equation} \label{phase}
\phi_1<\phi_2<\ldots<\phi_k<\phi_1+1.
\end{equation}

Given any two integers $i\le j$, we define $P_{ij}:=P_i\#P_{i+1}\#\ldots\#P_j$ as
in the Proposition above, and let $\phi_{ij}$ be the unique real number in the
interval $[\phi_i,\phi_j]\subseteq [\phi_1,\phi_1+1)$ such that
$m_ie^{i\pi\phi_i}+\ldots+m_je^{i\pi\phi_j}=m_{ij}e^{i\pi\phi_{ij}}$
for some positive real numbers $m_{ij}$. Note this gives us inequalities
\begin{equation} \label{phase2}
\phi_{ij}<\phi_{kl} \quad\text{if}\ i<k\ \text{and}\ j<l,\ \text{and}\
\phi_{ij}\in[\phi_1,\phi_1+1)\ \ \forall ij.
\end{equation}

\begin{defn} \label{stable}
Fix $N\ge2$ and $m_{ij},\,\phi_{ij}$ as above. Define
$$
Z(P_i):=m_ie^{i\pi\phi_i},
$$
and extend $Z$ to be defined on all of $K(D_k)=\oplus_i\,\Z_{P_i}$ by linearity,
so that $Z(P_{ij})=Z(P_i)+\ldots+Z(P_j)=m_{ij}e^{i\pi\phi_{ij}}$.

Then define the $P_{ij}[m]$s to be the stable objects of phase $\phi_{ij}+m$
in $D^N_k$, and define the $\SS(\phi)$s to
consist of all direct sums of the stable objects of phase $\phi$.
\end{defn}

To show this is a stability condition we need to understand the Homs between
these $P_{ij}$s, which we do now using Proposition \ref{homs}.

\begin{prop} \label{morehoms}
For $i<k<j+1<l+1$ we have
\begin{equation} \label{ijkl}
\Hom(P_{kl},P_{ij})\cong\C\,[-1]\oplus\C\,[-N],
\end{equation}
while if one of the inequalities becomes an equality we get only one of the
two summands: for $i=k<j+1<l+1$ and $i<k<j+1=l+1$ we get $\C\,[-N]$, while
for $i<k=j+1<l+1$ we get $\C\,[-1]$.

The duality $\Hom(P_{kl},P_{ij})\cong\Hom(P_{ij},P_{kl})^\vee[-N]$
(\ref{dual}) applied to (\ref{ijkl}) determines more Homs. All others are zero
apart from $\Hom(P_{ij},P_{ij})\cong\C\oplus\C\,[-N]$:
the $P_{ij}$ are spherical.
\end{prop}

\begin{pf} Building up $P_{ij}=P_{i,j-1}\#P_j=P_i\#P_{i+1,j}$ inductively
and using Proposition \ref{homs} gives this result very easily. We give
the example of most interest to us: $i<k=j+1<l+1$.

Using $\Hom(P_r,P_{ij})\cong0$ for $r>j+1$ (\ref{homs}), it is easy to
show inductively that $\Hom(P_{rl},P_{ij})=0$ for $r>j+1$ (where if $l<r$
we define $P_{rl}:=0$). Then applying $\Hom(\ \cdot\ ,P_{ij})$ to
$$
P_{j+1}\to P_{j+1,l}\to P_{j+2,l}
$$
gives (\ref{homs})
$$
\C\,[-1]\to\Hom(P_{j+1,l},P_{ij})\to0\to\C.
$$
Thus $\Hom(P_{j+1,l},P_{ij})=\C\,[-1]$ as required.

The fact that the $P_{ij}$ are spherical is also proved inductively from the observation
that if $A$ and $B$ are spherical with $\Hom(B,A)=\C\,[-1]$, then the corresponding
extension $A\#B$ is also spherical (the connect
sum of two spheres is a sphere!): from the triangle $A\to A\#B\to B$
we get
$$
\xymatrix@R=.85em@C=1em{
\Hom(B,A) \ar[r] \ar[d] & \Hom(B,A\#B) \ar[r] \ar[d] & \Hom(B,B) \ar[d] \\
\Hom(A\#B,A) \ar[r] \ar[d] & \Hom(A\#B,A\#B) \ar[r] \ar[d] & \Hom(A\#B,B) \ar[d] \\
\Hom(A,A) \ar[r] & \Hom(A,A\#B) \ar[r] & \Hom(A,B).
}$$
In the first column the connecting map $\Hom^0(A,A)\to\Hom^1(B,A)$ takes the
identity to the generator, by definition of the nontrivial extension. Since
$A$ is spherical, this makes $\Hom(A\#B,A)\cong\C\,[-N]$. Similarly with
the last column, using the functorial duality $\Hom(E,F)^\vee\cong\Hom(F,E)[N]$.
So the central row becomes
$$
\C\,[-N]\to\Hom(A\#B,A\#B)\to\C,
$$
from which it follows that $A\#B$ is spherical, since $N\ge2$.
\end{pf}

Finally we can prove that we are in the situation of Theorem \ref{split}.

\begin{thm}
Definition \ref{stable} defines a stability condition (\ref{def}) on $D^N_k$.
\end{thm}

\begin{figure}[h]
\center{
\input{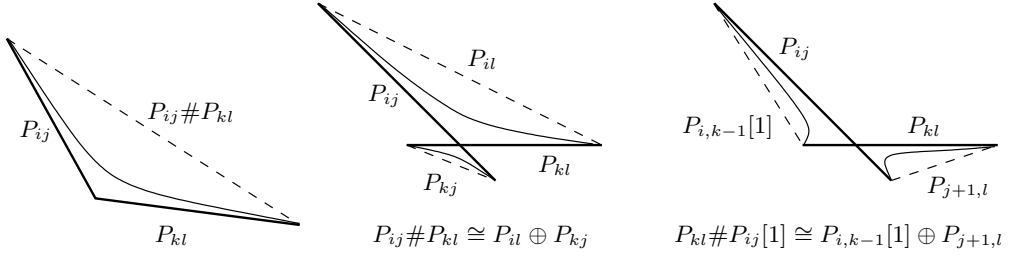}
\caption{Connect sums of stable objects $P_{ij}$ and $P_{kl}$ of increasing phase
are (direct sums of) stable objects: they (the curved lines) can be deformed to the
dashed straight lines. \label{connect}}}
\end{figure}

\begin{pf} 
Axioms (a), (b) and (e) of (\ref{def}) are immediate from Definition \ref{stable},
while (d) follows from Proposition \ref{morehoms} and the inequalities
(\ref{phase2}).
Axiom (c) will follow from Theorem \ref{split} if we can show that any
\emph{nontrivial} extension $E\to C\to F$, with $E,\,F$ stable of phases
$\phi<\psi$, $C$ is either stable (of phase $\theta\in(\phi,\psi)$),
or a sum of stables $C=A\oplus B$, with
$A\in\S(\alpha),\,B\in\S(\beta),\ \phi<\alpha\le\beta<\psi$.

Such extensions are given by the shifts of the Homs computed in Proposition
\ref{morehoms}. Only the Homs listed there in degrees 0 and 1 interest
us: if there is a nonzero $\Hom^n(F,E)$ with $n\ge2$ and the phases of
$E$ and $F$ less than 1 apart (\ref{phase2}) then this gives rise to an extension
in $\Hom^1(F,E[n-1])$ with the phase of $E[n-1]$ \emph{greater} than the phase
of $F$, which therefore does not concern us. Note, however, that by duality,
Homs of degree $N,\,N-1$ give rise to Homs of degree $0,\,1$, in the opposite
direction, that we do need to consider. This gives us 5 cases to check
(there are also selfHoms to consider, but these just give $P\#P[1]\cong0$).

\noindent $\bullet$ $P_{ij}\#P_{j+1,l}$, and their shifts. We have already
observed that this gives the single stable object $P_{il}$ with phase
$\phi_{il}\in(\phi_{ij},\phi_{j+1,l})$.

\noindent $\bullet$ $P_{kl}\#(P_{kj}[1]),\ j<l,$ and their shifts. This extension
comes from the element of $\Hom^0(P_{kj},P_{kl})$ in the triangle
$P_{kj}\to P_{kl}\to P_{j+1,l}$, and so is isomorphic to $P_{j+1,l}$, with
phase $\phi_{j+1,l}\in(\phi_{kl},\phi_{kj}+1)$.

\noindent $\bullet$ $P_{kl}\#(P_{il}[1]),\ i<k,$ and their shifts. This extension
comes from the element of $\Hom^0(P_{il},P_{kl})$ in the triangle
$P_{i,k-1}\to P_{il}\to P_{kl}$, and so is isomorphic to $P_{i,k-1}[1]$, with
phase $\phi_{i,k-1}+1\in(\phi_{kl},\phi_{il}+1)$.

\noindent $\bullet$ $P_{ij}\#P_{kl},\ i<k<j+1<l+1,$ and their shifts. This
element of $\Hom^1$ is in the image of $\Hom^1(P_{kl},P_{i,k-1})$ (that
defines $P_{i,k-1}\#P_{kl}=P_{il}$) via the map $P_{i,k-1}\to P_{ij}$ (with
cone $P_{kj}$). So we induce a diagram
$$
\xymatrix@R=.85em{
P_{i,k-1} \ar[r] \ar[d] & P_{il} \ar[r] \ar[d] & P_{kl} \ar@{=}[d] \\
P_{ij} \ar[r] \ar[d] & P_{ij}\#P_{kl} \ar[r] \ar[d] & P_{kl} \\
P_{kj} \ar@{=}[r] & P_{kj}
}$$
where the bottom row is induced from the octahedral Lemma. So the central
column shows that $P_{ij}\#P_{kl}\cong P_{kj}\oplus P_{il}$, as there are
no nontrivial extensions between these two objects, by Proposition
(\ref{morehoms}). (See also Figure \ref{connect}.) The inequalities
(\ref{phase2}) give the required $\phi_{kj},\phi_{il}\in
(\phi_{ij},\phi_{kl})$.

\noindent $\bullet$ $P_{kl}\#(P_{ij}[1]),\ i<k<j+1<l+1,$ and their shifts.
Similarly the diagram
$$
\xymatrix@R=.85em@C=.85em{
P_{ij} \ar@{=}[r] \ar[d] & P_{ij} \ar[d] \\
P_{kj} \ar[r] \ar[d] & P_{kl} \ar[r] \ar[d] & P_{j+1,l} \ar@{=}[d] \\
P_{i,k-1}[1] \ar[r] & P_{kl}\#(P_{ij}[1]) \ar[r] & P_{j+1,l}
}$$
of the standard morphisms we have already seen shows that $P_{kl}\#P_{ij}[1]\cong
P_{j+1,l}\oplus P_{i,k-1}[1]$ since there are no extensions between these
two objects. Again (\ref{phase2}) shows that $\phi_{j+1,l},\phi_{i,k-1}+1\in
(\phi_{kl},\phi_{ij}+1)$.
\end{pf}

\section{Deforming the stability condition}

In this section we determine an entire connected component of the space
of stability conditions in dimension $N=2$ and connect it with braid groups
of autoequivalences \cite{RZ}, \cite{KS}, \cite{ST}. The proofs are a little brief; in particular they make use of the description of stability conditions
mentioned in (\ref{heart}) and proved in \cite{Br1}.

In \cite{Br1} the space of stability conditions is shown to be a metric space
in a natural way, such that the mass function (\ref{mass}) is continuous.
Using the description of stability described in the
paragraph (\ref{heart}), the space is shown locally, about any given
stability condition, to be isomorphic to the space $\Hom(K(\T),\C)$ of
$Z$s. In fact it is shown to be a cover of the space of $Z$s minus those
where a mass (\ref{mass}) goes to zero. That is the stability condition deforms
even through walls of $Z$s where some objects become unstable and others stable,
so long as $Z$ of a stable object does not become zero. We will show that
in our case this means that the points 0 and $Z(P_{1j})=\sum_{i=1}^jZ(P_i)$
are distinct, i.e., geometrically, the ``endpoints" of our special Lagrangians
are distinct -- no stable object's mass has gone to zero. On looping round
such a zero (a generator of the braid group) back to the same $Z$, we will
find that the stability condition has changed, as the set of
stable objects undergoes a ``Dehn twist" \cite{ST}, at least in the easiest
case to analyse, dimension $N=2$, to which we will restrict from now on.

To do this we will need show that under deformations of $Z$, our set of
stable objects keep the following properties of our initial set of stable
objects (\ref{stable}).

\begin{defn} \label{simple}
The set of stable objects of a stability condition on $D^2_k$ is called
\emph{simple} if it is the set of shifts of $k(k+1)/2$ distinct spherical
objects $Q_{ij},\ 1\le i\le j\le k$, satisfying the following conditions.
$[Q_{ij}]=[P_{ij}]$ in K-theory, there is a single Hom in some degree
between $Q_{ab}$ and $Q_{b+1,c}$, a single Hom in some degree
between $Q_{ab}$ and $Q_{ac}\ (c\ne a)$, and the Euler characteristic
\begin{equation} \label{euler}
\chi(Q_{ab},Q_{cd}):=\sum_i(-1)^i\dim\Hom^i(Q_{ab},Q_{cd})=0
\quad\text{for}\quad a<c\ne b+1.
\end{equation}
We call the
stability condition simple if its set of stable objects is simple and its
semistable objects are direct sums of stable objects of the same phase.
\end{defn}

We want to show next that the stability condition remains simple on
passing through walls of semistability.

We say that $Z\in\Hom(K(D_k),\C)$ lies on a \emph{codimension one wall}
if $Z([P_{ij}])\ne0\ \forall i,j$, at least 2 of the classes $[P_{ij}]$
have the same phase (mod 1), and at most 3, in which case the three must be linearly
dependent (and so of the form $[P_{ab}],\,[P_{b+1\,c}],\,[P_{ac}]$ -- notice
that for the sum of two $[P_{ij}]$ classes to equal the class of another,
the two classes must be of the form $[P_{ab}],\,[P_{b+1\,c}]$.)

We can talk locally about ``sides" of such a wall depending on the sign of
the difference in sign of these two phases.

\begin{prop} \label{onehom}
If $P,\,Q$ are distinct stable objects in a simple stability condition, close
to, and on one side of, a codimension one wall, and whose phases coincide
on the wall, then the total dimension of $\Hom^*(P,Q)$ is at most 1.

If the total dimension of $\Hom^*(P,Q)$ is exactly 1, then there is a further
stable object $R$ of the same phase on the wall, and, on reordering
$P,Q,R$ if necessary, the Hom is of degree 1 and $R\cong P\#Q$. 
\end{prop}

\begin{pf} We need to rule out there
being two or more Homs between $P$ and $Q$. In this case, by the definition of
codimension one wall and simple (\ref{simple}), there are no other K-theory
classes of stable objects whose phases tend to those of $[P],\,[Q]$ on
the wall.

Shifting $Q$ and swapping $P,\,Q$ if necessary, and moving closer to the
wall, we can assume that the
phase of $Q$ is more than that of $P$, and that there are no stable objects
of phase in between. Thus, by stability
and Serre duality, $\Hom^i(Q,P)=0$ for $i\le0$ and $i\ge3$. By the vanishing
of the Euler characteristic (\ref{euler}), then, there must be equal numbers
of Homs in degrees 1 and 2; pick a nonzero element of $\Hom^1(Q,P)$ and
form the extension $P\#Q$. Using the description of stability given in
the paragraph (\ref{heart}),
if $P\#Q$ were unstable there would be a stable object, with a nonzero Hom
to $P\#Q$, of phase between those of $P,\,Q$. But this is a contradiction,
so that $P\#Q$ is stable -- another contradiction, since its K-theory class
is $[P]+[Q]$, which does not contain a stable object.

So we need only consider the case where there is a single Hom from $Q$
to $P$, which by axiom (d) and Serre duality must be in degree 1 or 2 (for
the phase of $Q$ greater than that of $P$).
There are now 3 stable objects $P,Q,R$ whose phase is
tending to the same value on the wall. We work with the semistability
criterion in the abelian category (\ref{heart})
(for a suitable value of $t$ so that it contains $P,Q,R$). Without loss
of generality we will assume that the masses (\ref{mass}) of $P,Q$ are
less than or equal to that of $R$ on the wall.

There can be no element of $\Hom^2(Q,P)\cong\Hom^0(P,Q)^*$,
as the image in $Q$ of any Hom from $P$ (i.e.\ the cokernel of its kernel
in the abelian heart) would have a filtration by stable objects
of phase between those of $P$ and $Q$, and of strictly smaller mass.

So $\Hom^1(Q,P)=\C$, and we can form $P\#Q$. Its Harder-Narasimhan filtration
in the abelian category (\ref{heart}) is of semistable objects of phase
between those of $P$ and $Q$; i.e.\ of direct sums of shifts of $P,Q,R$.
Since the extension is non-trivial and the masses of $P,Q$ are less than
or equal to that of $R$, it follows that the filtration is just
$R$, that is $P\#Q=R$. We claim that there is a single Hom between $R$ and
either of $P$ and $Q$; for instance $\Hom(R,P)$ fits into long exact sequence
$$
\Hom^i(Q,P)\to\Hom^i(P\#Q,P)\to\Hom^i(P,P)\to\Hom^i(Q,P)\to\ldots,
$$
in which the identity in $\Hom(P,P)$ maps to the generator of $\Hom^1(Q,P)$.
Thus $\Hom(P\#Q,P)\cong\C\,[2]$. Similarly $\Hom(Q,P\#Q)\cong\C\,[2]$.
\end{pf}

\begin{prop}
On crossing a codimension one wall a simple stability condition remains
simple.
\end{prop}

\begin{pf}
Suppose the phases of two stable objects $P$ and $Q$
(on one side of the wall,
where the phase of $P$ is less than that of $Q$, without loss of generality)
coincide on the wall. If there are no Homs between them it is easy to see
that there are no other such stable objects, and the set of (semi)stable objects
does not change across the wall; in particular $P$ and $Q$ remain stable
(see for example \cite{Br1} for more details, or use the equivalent definition
of stability in the paragraph (\ref{heart})).

So by Proposition \ref{onehom} we need only consider the case where there
is a single Hom from $Q$ to $P$, of degree one, and $R=P\#Q$ is also stable.

Then for $P,Q$ to become unstable as we cross the wall, they must have
a filtration by stable objects of strictly smaller mass and the same phase
on the wall. But because this is a codimension one wall and $R$ has greater
mass, no such objects exist.

$R$ becomes unstable on the other side of
the wall. But there we can use Serre duality to form a unique
nontrivial extension $Q\#P$, which is similarly stable on that side of
the wall as in the proof of Proposition \ref{onehom}. $Q\#P$ is also spherical
as in the proof of Proposition \ref{morehoms}. Since no other
stable objects are affected we claim that this new set of stable
objects (with each shift $P\#Q[r]$ replaced by $Q\#P[r]$) is also simple.
The K-theory class of $Q\#P$ is the same as that of $P\#Q$, and its Homs
to $P$ and $Q$ are one dimensional as in the last paragraph of the proof
of Proposition \ref{onehom}. To satisfy Definition \ref{simple}, then,
we must finally check that for any stable object $E$ with $\chi(E,P\#Q)=0$
(\ref{euler}), we also have $\chi(E,Q\#P)=0$. But by the exact sequence
$\Hom^i(E,Q)\to\Hom^i(E,Q\#P)\to\Hom^i(E,P)\to\Hom^{i+1}(E,Q)\to\ldots$
we see that
$$
\chi(E,Q\#P)=\chi(E,Q)+\chi(E,P)=\chi(E,P\#Q)=0,
$$
as required.
\end{pf}

While we have analysed crossing only codimension one walls, the fact from
\cite{Br1}, mentioned above, that locally the space of stability conditions
is isomorphic to the space of $Z$s means that there is no monodromy around
codimension 2 walls, and whenever $Z$ does not lie on a wall, the stability
condition is \emph{simple}. In particular, then, there is always a stable
object in the K-theory class $[P_{ij}]$ away from the finite number of
walls. Combined with the results on deforming
stability conditions \cite{Br1}, which can always be done until a stable
object's mass goes to zero, we find the connected component of our stability
conditions (\ref{stable}) is a cover of the space of ${Z}$s such that
$Z([P_{ij}])\ne0$ for all $i\le j$. Plotting the points $0$ and $Z([P_{1i}])$
for all $i$, and translating them to have mean zero, this space in turn covers
the configuration space $C_{k+1}^0$ of $(k+1)$ \emph{distinct} points in
$\C$ with centre of mass the origin. We now prove slightly more.

\begin{thm}
Via the above map, the connected component of the stability conditions
(\ref{stable}) is the \emph{universal} cover of the configuration space
$C_{k+1}^0$. The deck transformations are given
by the $B_k$ action of \cite{KS}, \cite{ST}.
\end{thm}

\begin{pf}
We want to check that the result of going round loops in configuration
space is the braid group action of \cite{KS}, \cite{ST}. It is sufficient
to check this for a suitable choice of generators; namely we pick a stability
condition as in (\ref{stable}) in which the mass $|Z(P_{ii})|$ of the
stable object $P_i$ is strictly smaller than that
of all other stable objects, and we move $Z$ by rotating
$Z(P_i)$ anticlockwise through $\pi$ radians while fixing $Z$
of stable objects of different endpoints. As we rotate $Z(P_i)$, the phases
of the stable objects $P$ with an endpoint $(i-1)$ or $i$ will also change (by no
more than $\epsilon=\sin^{-1}(|Z(P_i)|/|Z(P)|)$, by the triangle inequality); we
also assume that $|Z(P_{ii})|$ is so small that no two phases of classes
$[P_{kl}]\ne[P_i]$ coincide under this rotation; i.e.\ we take all the phases
of the classes $[P_{kl}]$ to be distinct, and $\epsilon$ to be smaller than
the smallest difference in such phases.

Note that for this choice of stability condition and stable object $P_i$, there
are no other stable objects with 2 Homs to $P_i$ (\ref{homs}); they either have
no Homs or one. As we rotate $Z(P_i)$ and follow the stability condition below,
we can see that this property is preserved, as by design no stable objects'
phases become equal, except to the phase of $P_i$, so no stable objects
change except via their interaction with $P_i$ described now.

As we rotate $Z(P_i)$ and cross walls (of codimension one only, without
loss of generality, by perturbing the loop if necessary), stable objects
with no Homs are unaffected, while those with one Hom are altered as in
the proof of the last Proposition.

That is, as the phase of $Z(P_i)$ reaches that of some other pair of stable objects
$E,F$ ($F$ of smaller mass) with $\Hom(P_i,E)\cong\C\,[0],\ \Hom(P_i,F)\cong\C\,[-1]$
and $E\cong P_i\#F$. On passing through the wall the stable object
$E$ is replaced by $F\#P_i$, and all other stable objects are left unchanged.

But we claim that
$$
T_{P_i}(P_i\#F)\cong F \qquad\text{and}\qquad T_{P_i}F\cong F\#P_i,
$$
where $T_{P_i}$ is the \emph{Dehn twist} about the spherical object $P_i$;
an equivalence of triangulated categories \cite{ST}. $T_{P_i}(E)$ sits in the
triangle
$$
P_i\otimes\Hom(P_i,E)\to E\to T_{P_i}(E),
$$
where the first map is evaluation. Since $\Hom(P_i,P_i\#F)\cong\C\,[0]$ and
$\Hom(P_i,F)\cong\C\,[-1]$, applying this to $E=P_i\#F$ and $E=F$
gives the standard triangles
$$
P_i\to P_i\#F\to F \qquad\text{and}\qquad P_i[-1]\to F\to F\#P_i
$$
respectively, proving the claim.

So on listing the stable objects
in ascending phase, we find that the subsequence $P_i,E,F$ is replaced
by $F,F\#P_i,P_i$, i.e.\ by $T_{P_i}E,T_{P_i}F,P_i$.

The same is true when the phase of $P_i$ passes through that of the other
stable objects $P$ with no Homs to $P_i$: they are left unaffected, just
as under $T_{P_i}$: $T_{P_i}P\cong P$.

As we rotate $Z(P_i)$ through $\pi$, its phase crosses the phase of \emph{all}
stable objects (after a suitable shift), and so we end up with the list
of stable objects of ascending phase $P_i,A,B,\ldots$ being replaced by
$T_{P_i}A,T_{P_i}B,\ldots,
P_i$. Equivalently, as $P_i$s phase has increased by $1$, we have $P_i[-1],
T_{P_i}A,T_{P_i}B,\ldots$. But $T_{P_i}P_i\cong P_i[-1]$, so we have altered
the set of stable objects (and $Z$) by the action of $T_{P_i}$, as claimed.

Finally, since the $B_k$ action of (\cite{KS}, \cite{ST}) is faithful not
just on the triangulated categories, but also on $A_k$-chains such as the
stable objects $\{P_i\}$ (\cite{ST} Theorem 4.13), we see that the cover
of $C_{k+1}$ that we get is the \emph{universal} cover.
\end{pf}

We end by noting that the if we take an element of the braid group which
acts trivially on K-theory (and so on $Z$), then the two stability conditions
differ only by their set of stable objects. In our case, these stable objects
differed by an autoequivalence of the triangulated category. More generally,
if the axiomatic notion of stability of \cite{Br1} is to agree with the
physical notion of stability (and the geometric conjectures in \cite{Th})
this would have to hold more generally; that is one might conjecture that
two stability conditions with the same central charge on a ``Calabi-Yau
category" (one with a ``trace map" $\Hom^N(E,E)\to\C$ inducing a functorial
duality $\Hom(A,B)^*\cong\Hom(B,A)[N]$ for all $A,B$) should differ by an
autoequivalence of the category that is the identity on its numerical K-theory.

\vspace{4mm}
\begin{flushleft}
Department of Mathematics, Imperial College, 180 Queen's Gate,
London SW7 2BZ. UK. \\
{\small Email: {\tt rpwt@ic.ac.uk}}
\end{flushleft}

\end{document}